# A Quasi Time Optimal Receding Horizon Control

Piotr Bania
AGH University of Science and Technology,
Faculty of Electrical Engineering, Automatics, Informatics and Electronics
Department of Automatics.
Mickiewicza 30 av. 30-059 Krakow, Poland, fax: (+4812) 6341568
e-mail: pba@ia.agh.edu.pl

**Abstract:** This paper presents a quasi time optimal receding horizon control algorithm. The proposed algorithm generates near time optimal control when the state of the system is far from the target. When the state attains a certain neighbourhood of the aim, it begins the adaptation of the cost function. The purpose of this adaptation is to move from the time optimal control to the stabilizing control. Sufficient conditions for the stability of the closed loop system and the manner of the adaptation of the cost function have been given. Considerations are illustrated with examples.

**Introduction.** Observations of purposeful activities of animals, men and various biological systems lead to the conclusion that in these systems control aim is attained in two stages. In the first stage the undertaken decisions are to guarantee a quick, time optimal attainment of certain neighbourhood of the target. In the second phase, as the target is nearing, the control becomes more and more precise (soft, economical). At this stage accuracy is the priority. Two different quality criteria correspond to the two phases. In the first phase, the time to reach the target is essential. In the second stage, the minimization of state and control deviations is of the utmost importance. The behaviour described above can be also substantiated through analysis of the closed loop system sensitivity to disturbances and model uncertainty. It is known [12] that time optimal control is very sensitive to disturbances and modelling errors. The greatest sensitivity usually appears in the latter part of the process, when the trajectory overlies the switching surface [12]. It is also known that linear-quadratic algorithms (LQR linear quadratic regulator) and their non-linear versions (NLQR nonlinear quadratic regulator) are a lot less sensitive to disturbances and modelling errors. This decrease of sensitivity usually takes place at the expense of the control time. Therefore, it can be expected that the transition made in an adequately regular manner from the time optimal control to the stabilizing control will enable us to reach the compromise between speed, precision and sensitivity. In former papers the authors focused on algorithms in which the time horizon and the cost function were fixed [2, 3, 4, 5, 6]. In "dual–mode approach" proposed by Michalska [10], the control action is generated by linear feedback controller, when state of the system belongs to the neighbourhood of the origin. Outside this neighbourhood receding horizon control is employed and variable time horizon is permitted. This paper expands these results and proposes a Quasi Time Optimal Receding Horizon Control algorithm (QTO–RHC). In this approach cost function consists of three terms. The first term represents time horizon. The time horizon is a decisive variable. When the state of the system reaches neighbourhood of the origin, the algorithm gradually 'inserts' a second integral term into the cost function. The third term is so called terminal penalty term. The paper consists of six chapters. The first chapter gives basic definitions and theorems. The optimal control problem and the description of the general algorithm comprise the second chapter. The third chapter presents the proof of the closed loop system stability and gives certain auxiliary theorems which will permit to design a QTO–RHC algorithm that combines advantages of the time optimal and stabilizing control. Properties of the QTO–RHC are described in the fourth chapter. The fifth chapter analyses the case of the quadratic cost function for systems which are stabilizable in the neighbourhood of zero [2, 1]. Examples of the control of nonlinear systems are shown in the sixth chapter.

**Notation.** The dot product and the norm of the vector in $R^n$ are denoted as $\langle x, y \rangle_H = x^T H y$, $\| x \|_H = \sqrt{x^T H x}$, where $H = H^T > 0$. The set of non-negative real numbers is marked as $R_0^+$. The expression $PC([0,T], R^m)$ (or $PC([0,\infty), R^m)$) denotes the space of the $R^m$ valued piecewise continuous functions in $[0,T]$ (or $[0,\infty)$) with norm $\| u \|_\infty = \operatorname*{esssup}_{t \geq 0} \| u(t) \|$.

**1. Equations of the system and the properties of solutions.** Assume that the controlled plant is described by a system of ordinary differential equations

$$\dot{x}(t) = f(x(t), u(t)), \quad x(0) = x_0, \quad t \in R_0^+. \tag{1.1}$$

The function $f: R^n \times R^m \to R^n$ is of $C^1$ class with respect to both arguments and $f(0,0) = 0$. Moreover, $f$ fulfils the Lipschitz condition $\| f(x_1, u) - f(x_2, u) \| \leq L \| x_1 - x_2 \|$ in $R^n$, with constant $L > 0$ independent of the choice of points $x_1, x_2, u$. Constraints of the control values have the following form

$$u(t) \in U, \quad U = \{ u \in R^m; u_{\min} \leq u \leq u_{\max}, u_{\min} < 0, u_{\max} > 0 \} \tag{1.2}$$

The state of the system (1.1) in moment $t$ will be also denoted by the symbol $x_t$ i.e. $x_t = x(t)$. Let $\delta > 0$ represent the sampling time. In moments $t_i = i\delta$, $i = 0,1,2,..$ we estimate the state of the plant and we set control. The time horizon and the minimal time horizon are represented by $T$ and $T_{\min}$. The time horizon is a decisive variable, while the minimal horizon is fixed. Let's assume that both horizons satisfy

$$T \geq T_{\min} \geq \delta > 0 \tag{1.3}$$

The trajectory in the time interval $[t_i, t_i + T]$ is calculated by solving equations (1.1) with the initial condition $x(t_i) = x_{t_i}$. Solution of the equation (1.1) with control $u$, starting from point $x_t$ will be marked as $x(s; x_t; u(\cdot))$. In certain cases instead of $x(s; x_t; u(\cdot))$ we will write shortly $x(s)$. Let $\Omega \subset R^n$ be a closed, bounded and simply connected point set and let $0 \in \Omega$. Let us define admissible control.

**Definition 1.1.** Let $x$ be the solution of the equation (1.1) with control $u: [t, t+T] \to R^m$ and initial condition $x_t$. The control $u$ is admissible in point $x_t$ if: **1.** $u \in PC([t, t+T], R^m)$, **2.** $u(s) \in U$, $s \in [t, t+T]$, **3.** $x(t+T) \in \Omega$. □

The set of all admissible controls in point $x_t$ is marked as $\mathbf{U}(x_t)$.

**Definition 1.2.** If for every $x_0 \in R^n$ and for any number $\varepsilon > 0$ there exists an admissible control $u^D: [0, T_D] \to U$, $T_D < +\infty$ such that $\| x(T_D; x_0; u^D(\cdot)) \| < \varepsilon$, then the system (1.1) will be called asymptotically controllable. □

We assume that the system (1.1) is asymptotically controllable.

**Definition 1.3.** If the function $\varphi: R_0^+ \to R_0^+$ is continuous, strictly increasing and $\varphi(0) = 0$, then we speak of this function as class $\mathcal{K}$ function or we write $\varphi \in \mathcal{K}$. □

**Theorem 1.1.** Let $u: [t, t+T] \to R^m$ be the admissible control. With above assumptions holding, then the solution $x(s; x_t; u(\cdot))$ of equation (1.1) fulfils estimation

$$\| x(s) \| \leq (M_1 + M_2 T) e^{LT} \text{ where } M_1 = \sup_{z_0 \in \Omega} \| z_0 \|, \quad M_2 = L \sup_{z_0 \in \Omega} \| z_0 \| + \sup_{(z_0, w) \in \Omega \times U} \| f(z_0, w) \|. \square$$

For proof of the theorem see [14].

**2. Optimal control problem.** Controller's activity consists in a cyclical solution of the optimal control problem. In order to calculate control in time interval $[t_i, t_i + \delta]$, we solve the following problem.

**Problem** $P(t_i, x_{t_i})$: To find admissible control $u \in \mathbf{U}(x_{t_i})$, and final time $T \geq T_{\min}$, that minimize the cost function

$$J(u,T;x_{t_i}) = T + \varepsilon_i \int_{t_i}^{t_i+T} L(x(s),u(s))ds + \rho_i q(x(t_i+T)), \qquad (2.1)$$

where $x$ fulfils equation (1.1) with control $u$ and initial condition $x_{t_i}$. □

We assume that $L \in C^1(R^n \times R^m, R_0^+)$, $L(0,0) = 0$, and that $L(x,u) \geq \alpha_L(\|x\|)$ where $\alpha_L \in \mathcal{K}$. The function $q \in C^1(R^n, R_0^+)$, $q(0) = 0$ and there exists such a constant $c > 0$ that $q(x) \geq c\|x\|_H^2$. The sequences $\varepsilon_i$, $\rho_i$ fulfil the following conditions:

$$\varepsilon_i \in [0,1], \quad \rho_i \geq 1. \square \qquad (2.2)$$

The solution of the problem $P(t_i, x_{t_i})$ and the corresponding trajectory are marked as $\bar{u}^i(s, x_{t_i})$, $\bar{T}_i$, $\bar{x}^i(s; x_{t_i}; \bar{u}^i(s, x_{t_i}))$. We assume that solution of the above problem exists in the entire space or in an adequately spacious set $X$ that includes the set $\Omega$. The calculations are repeated every time $\delta$ for a new, currently given state of the plant. In certain cases we will skip arguments and write e.g. $\bar{u}^i$, $\bar{x}^i$, instead of $\bar{u}^i(s, x_{t_i})$, $\bar{x}^i(s; x_{t_i}; \bar{u}^i(s, x_{t_i}))$. We will now show the schema for receding horizon control.

**Schema 2.1.**
0. Substitute $i = 0$, $t_i = 0$;
1. Make estimation of the state $x_{t_i}$;
2. Calculate numbers $\varepsilon_i$, $\rho_i$;
3. Solve $P(t_i, x_{t_i})$ and find $\bar{u}^i : [t_i, t_i + \bar{T}_i] \to U$;
4. Apply the initial part of the control $\bar{u}^i$ in time interval $t \in [t_i, t_i + \delta)$;
5. Substitute $t_{i+1} \leftarrow t_i + \delta$, $i \leftarrow i+1$ and go to 1.

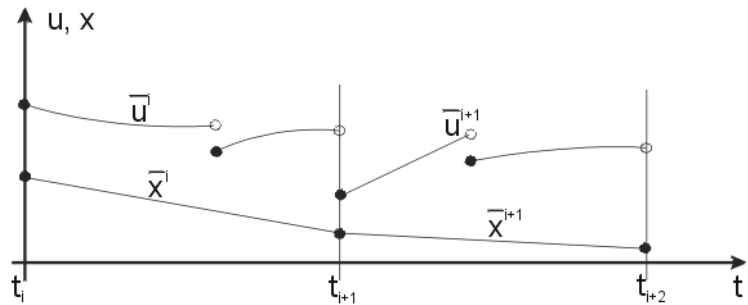

Fig. 2.1 Time structure in schema 2.1.

We obtain trajectory and control in the closed loop system through the concatenation of trajectory $\bar{x}^i$ and control $\bar{u}^i$ applied to the plant in time intervals $[t_i, t_{i+1})$, $i = 0, 1, 2, \ldots$ This trajectory and its corresponding control will be marked $x^*(t)$, $u^*(t)$. The trajectory $x^*(t)$ is continuous, while the control $x^*(t)$ is piecewise continuous.

**Conclusion 2.1.** Trajectory $x^*(t)$ has an infinite escape time.

*Proof:* From theorem 1.1 we have $\|\bar{x}^i(t)\| \leq (M_1 + M_2 \bar{T}_i) \exp(L \bar{T}_i)$ for $t \in [t_i, t_i + \bar{T}_i)$. Hence $\|x^*(t)\| \leq (M_1 + M_2 \bar{T}_i) \exp(L \bar{T}_i)$ for $t \in [t_i, t_i + \delta)$. □.

**3. Stability of the schema 2.1.** Let us analyse now stability of schema 2.1. Certain ideas in the argumentation of the stability are borrowed from the works of Fontes, Findeisen and Chen ([6], [3], [2]).

**Definition 3.1.** The closed loop system will be called stable if for every $x_0 \in X$ (or $x_0 \in R^n$) the condition $\lim_{t \to \infty} x^*(t) = 0$ is fulfilled. □

**Definition 3.2.** Let $\bar{u}, \bar{T}$ be solution of $P(t, x_t)$. We call the optimal value of the cost function (2.1) for the initial condition $x_t$ i.e. $V(x_t) = J(\bar{u}, \bar{T}; x_t)$ the value function $V(x_t)$. □

We will give now sufficient conditions for the closed loop system stability. Let $\bar{u}^i$, $\bar{T}_i$ be solution of problem $P(t_i, x_{t_i})$ for certain values $\varepsilon_i$, $\rho_i$. Increases of sequences $\varepsilon_i$, $\rho_i$ will be marked $\Delta \varepsilon_i = \varepsilon_{i+1} - \varepsilon_i$, $\Delta \rho_i = \rho_{i+1} - \rho_i$.

**Theorem 3.1**. Assume all former conditions hold and:

1. Set $\Omega$ and function $q$ will be so chosen that for every $x_s(0) \in \Omega$ there exists a piecewise continuous control $u_s : [0, \delta] \to U$ that trajectory $x_s$ of the system (1.1) generated by this control remains in set $\Omega$ and satisfies the condition

$$\frac{d}{d\tau} q(x_s(\tau)) + L(x_s(\tau), u_s(\tau)) \leq 0, \; \tau \in [0, \delta]. \tag{3.1}$$

2. Problem $P(0, x_0)$ has a solution.

3. Sequences $\varepsilon_i$, $\rho_i$ fulfil inequalities

$$\Delta \varepsilon_i \int_{t_i+\delta}^{t_i+\bar{T}_i} L(\bar{x}^i, \bar{u}^i) ds + \Delta \rho_i q(\bar{x}^i(t_i + \bar{T}_i)) \leq \delta + (1-\xi)\varepsilon_i \int_{t_i}^{t_i+\delta} L(\bar{x}^i, \bar{u}^i) ds \text{ for } \bar{T}_i \geq T_{\min} + \delta \tag{3.2}$$

and

$$\Delta \varepsilon_i \left( \int_{t_i+\delta}^{t_i+\bar{T}_i} L(\bar{x}^i, \bar{u}^i) ds + \int_{t_i+\bar{T}_i}^{t_i+T_{\min}+\delta} L(x_s, u_s) ds \right) + \Delta \rho_i q(x_s(t_i + T_{\min} + \delta)) \leq (\bar{T}_i - T_{\min}) +$$

$$+ (1-\xi)\varepsilon_i \int_{t_i}^{t_i+\delta} L(\bar{x}^i, \bar{u}^i) ds, \tag{3.3}$$

for $\bar{T}_i \in [T_{\min}, T_{\min} + \delta)$, but the number $\xi \in (0,1)$ is fixed.

4. There exists numbers $i_0 \geq 0$, $\varepsilon_{\min} > 0$ such that $\varepsilon_{\min} \leq \varepsilon_i \leq 1$ for $i \geq i_0$.

Then trajectory $x^*(t)$ asymptotically converges to zero.

*Proof*: Assume that problem $P(t_i, x_{t_i})$ has a solution $\bar{u}^i$. Utilization of control $\bar{u}^i$ in the time interval $[t_i, t_i + \sigma]$, $\sigma \in (0, \delta]$ gives us $x^*(t_i + \sigma) = \bar{x}^i(t_i + \sigma) \in X$. We will show that control

$$\tilde{u}^i(t) = \begin{cases} \bar{u}^i(t), \; t \in [t_i + \sigma, t_i + \bar{T}_i], \\ u_s(t), \; t \in [t_i + \bar{T}_i, t_i + \bar{T}_i + \sigma], \end{cases}$$

fulfils the conditions of definition 1.1 and is admissible in point $\bar{x}_{t_i+\sigma}$. Control $\tilde{u}$ is a concatenation of piecewise continuous controls therefore it is a piecewise continuous control. Admissibility of $\tilde{u}$ results in $\bar{x}(t_i + \bar{T}_i) \in \Omega$. Assumption 1 of theorem 3.1 indicates that there exists such control $u_s(t)$, $t \in [t_i + \bar{T}_i, t_i + \bar{T}_i + \sigma]$ that $x_s(t) \in \Omega$, for $t \in [t_i + \bar{T}_i, t_i + \bar{T}_i + \sigma]$. Using the assumption 2 of the theorem and applying inductive reasoning will show us that the existence of a solution to problem $P(0, x_0)$ implies the existence of a solution to all problems $P(t, x_t)$ for $t \geq 0$. We will show now that the value function decreases from certain moment in time, and the closed loop trajectory converges to zero. We will adapt the following reasoning. If the pair $\bar{u}^i$, $\bar{T}_i$, is the solution of $P(t_i, x_{t_i})$ then for any admissible solution $\tilde{u}^i$, $\tilde{T}_i$ inequality $V(x_{t_i}) \leq J(\tilde{u}^i, \tilde{T}_i; x_{t_i})$ occurs. Control $\tilde{u}^i$ will be constructed in the way described above, or through limiting the control $\bar{u}^i$ to the interval $[t_i + \sigma, t_i + \overline{T}_i]$. Assumption 4 warrants that from certain moment $i_0\delta$ the condition $\varepsilon_{\min} \leq \varepsilon_i \leq 1$ is fulfilled. Conclusion 2.1 claims that trajectory $x^*$ is bounded in time interval $[0, i_0\delta]$, we can however say nothing about its convergence to zero. We will consider two cases.

**Case A.** Let $\varepsilon_{\min} \leq \varepsilon_i \leq 1$ and $\overline{T}_i \geq T_{\min} + \delta$, (Fig.3.1).

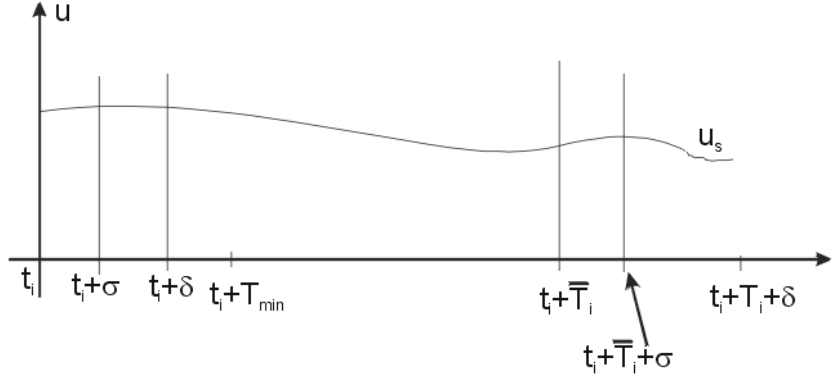

Fig.3.1. Time structure of control in case A.

The value function for the start from point $x^*_{t_i}$ equals

$$V(x^*_{t_i}) = \overline{T}_i + \varepsilon_i \int_{t_i}^{t_i+\overline{T}_i} L(\bar{x}^i, \bar{u}^i)ds + \rho_i q(\bar{x}^i(t_i + \overline{T}_i)). \tag{A.0}$$

We will show that $V(x^*_{t_i})$ decreases in interval $[t_i, t_i + \delta)$. Let $\bar{u}^\sigma$, $\overline{T}_\sigma$ represent a solution to $P(t_i + \sigma, x^*_{t_i+\sigma})$. Since the pair $\bar{u}^\sigma$, $\overline{T}_\sigma$ minimizes $J(u;T;x^*_{t_i+\sigma})$ hence for $\sigma \in (0,\delta)$ the following inequalities occur

$$V(x^*_{t_i+\sigma}) = \overline{T}_\sigma + \varepsilon_i \int_{t_i+\sigma}^{t_i+\overline{T}_\sigma+\sigma} L(\bar{x}^\sigma, \bar{u}^\sigma)ds + \rho_i q(\bar{x}^\sigma(t_i + \sigma + \overline{T}_\sigma)) \leq \overline{T}_i - \sigma + \varepsilon_i \int_{t_i+\sigma}^{t_i+\overline{T}_i} L(\bar{x}^i, \bar{u}^i)ds + \rho_i q(\bar{x}^i(t_i + \overline{T}_i))$$

and $V(x^*_{t_i+\sigma}) - V(x^*_{t_i}) \leq -\sigma - \varepsilon_i \int_{t_i}^{t_i+\sigma} L(\bar{x}^i, \bar{u}^i)ds \leq -\xi\varepsilon_i \int_{t_i}^{t_i+\sigma} \alpha_L(\|\bar{x}^i\|)ds < 0, \bar{x}^i \neq 0.$ (A.1)

After time $\delta$ we perform new optimization starting from point $x^*_{t_i+\delta}$. The value function equals now $V(x^*_{t_i+\delta}) = \overline{T}_{i+1} + \varepsilon_{i+1} \int_{t_i+\delta}^{t_i+\overline{T}_{i+1}+\delta} L(\bar{x}^{i+1}, \bar{u}^{i+1})ds + \rho_{i+1} q(\bar{x}^{i+1}(t_i + \overline{T}_{i+1} + \delta))$.

Since the pair $\bar{u}^{i+1}$, $\overline{T}_{i+1}$ minimizes $J(u;T;x^*_{t_i+\delta})$ then we obtain inequality

$$V(x^*_{t_i+\delta}) \leq \overline{T}_i - \delta + \varepsilon_{i+1} \int_{t_i+\delta}^{t_i+\overline{T}_i} L(\bar{x}^i, \bar{u}^i)ds + \rho_{i+1} q(\bar{x}^i(t_i + \overline{T}_i)).$$

Estimation of the difference in the value function results in

$$V(x^*_{t_i+\delta}) - V(x^*_{t_i}) \leq -\delta + \varepsilon_{i+1} \int_{t_i+\delta}^{t_i+\overline{T}_i} L(\bar{x}^i, \bar{u}^i)ds + \rho_{i+1} q(\bar{x}^i(t_i + \overline{T}_i)) - \varepsilon_i \int_{t_i}^{t_i+\overline{T}_i} L(\bar{x}^i, \bar{u}^i)ds - \rho_i q(\bar{x}^i(t_i + \overline{T}_i)) \leq$$

$\leq -\delta - \varepsilon_i \int_{t_i}^{t_i+\delta} L(\bar{x}^i, \bar{u}^i)ds + \Delta\varepsilon_i \int_{t_i+\delta}^{t_i+\overline{T}_i} L(\bar{x}^i, \bar{u}^i)ds + \Delta\rho_i q(\bar{x}^i(t_i + \overline{T}_i))$. From (3.2) we have

$$V(x^*_{t_i+\delta}) - V(x^*_{t_i}) \leq -\xi\varepsilon_i \int_{t_i}^{t_i+\delta} L(\bar{x}^i, \bar{u}^i)dt \leq -\xi\varepsilon_i \int_{t_i}^{t_i+\delta} \alpha_L(\|\bar{x}^i\|)dt < 0, \bar{x}^i \neq 0. \tag{A.2}$$

**Case B.** Let $\varepsilon_{\min} \leq \varepsilon_i \leq 1$ and $\overline{T}_i \in [T_{\min}, T_{\min} + \delta)$. In this case it is not possible to decrease the horizon of $\delta$. There are two possibilities.

**B.1.** Let $0 < \sigma \leq \overline{T}_i - T_{\min}$, (Fig.3.2).

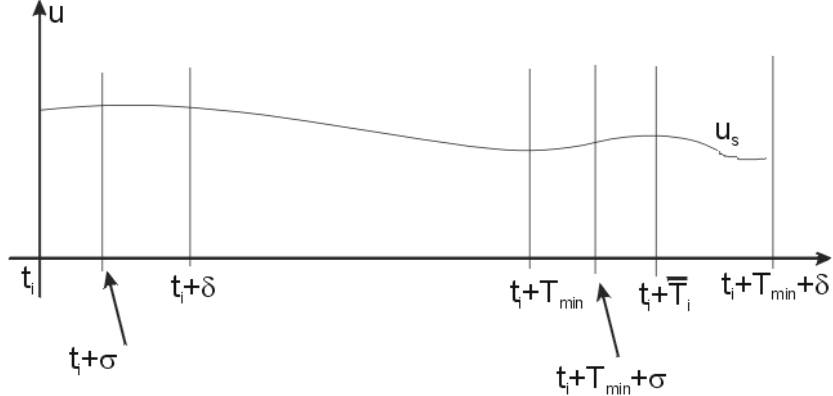

Fig. 3.2. Time structure of control in case **B.1**.

Adapting reasoning analogical to the one used in case **A** will give us again inequality **A.1**.

**B.2.** Let $0 \leq \overline{T}_i - T_{\min} < \sigma < \delta$ now, (see Fig.3.3) and let $u_s$ fulfil assumption 1. Consider the following control $\tilde{u}^i(t) = \begin{cases} \overline{u}^i(t), & t \in [t_i + \sigma, t_i + \overline{T}_i], \\ u_s(t), & t \in [t_i + \overline{T}_i, t_i + T_{\min} + \sigma]. \end{cases}$ This control is defined in the interval $[t_i + \sigma, t_i + T_{\min} + \sigma]$ and is admissible according to definition 1.1.

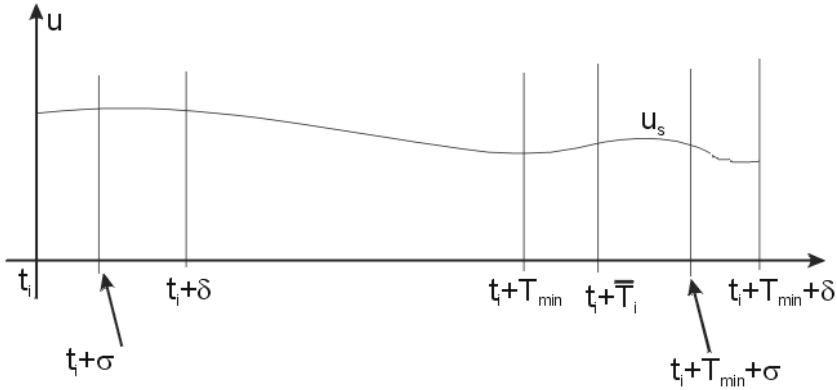

Fig.3.3. Time structure of control in case **B. 2**.

The following inequality occurs

$$V(x^*_{t_i+\sigma}) \leq T_{\min} + \varepsilon_i \int_{t_i+\sigma}^{t_i+\overline{T}_i} L(\overline{x}^i, \overline{u}^i) ds + \varepsilon_i \int_{t_i+\overline{T}_i}^{t_i+T_{\min}+\sigma} L(x_s, u_s) ds + \rho_i q(x_s(t_i + T_{\min} + \sigma)).$$

Let us estimate the difference in the value function

$$V(x^*_{t_i+\sigma}) - V(x^*_{t_i}) \leq -(\overline{T}_i - T_{\min}) - \varepsilon_i \int_{t_i}^{t_i+\sigma} L(\overline{x}^i, \overline{u}^i) ds + \varepsilon_i \int_{t_i+\overline{T}_i}^{t_i+T_{\min}+\sigma} L(x_s, u_s) ds +$$

$+ \rho_i \big(q(x_s(t_i + T_{\min} + \sigma)) - q(\overline{x}^i(t_i + \overline{T}_i))\big)$. From assumption 1 of the theorem and the fact that

$\varepsilon_i \in [0,1]$, $\rho_i \geq 1$ we obtain inequality $\varepsilon_i \int_{t_i+\overline{T}_i}^{t_i+T_{\min}+\sigma} L(x_s, u_s) ds + \rho_i \big(q(x_s(t_i + T_{\min} + \sigma)) - q(\overline{x}^i(t_i + \overline{T}_i))\big) \leq 0$.

Combining this two inequalities we have

$$V(x^*_{t_i+\sigma}) - V(x^*_{t_i}) \leq -(\overline{T}_i - T_{\min}) - \varepsilon_i \int_{t_i}^{t_i+\sigma} L(\overline{x}^i, \overline{u}^i) ds \leq -\xi \varepsilon_i \int_{t_i}^{t_i+\sigma} \alpha_L(\|\overline{x}^i\|) ds < 0, \; \overline{x}^i \neq 0.$$

In conclusion, for $\sigma \in (0,\delta)$ and $\overline{T}_i \in [T_{\min}, T_{\min} + \delta)$ the following estimation occurs

$$V(x^*_{t_i+\sigma}) - V(x^*_{t_i}) \leq -\min(\overline{T}_i - T_{\min}, \sigma) - \varepsilon_i \int_{t_i}^{t_i+\sigma} L(\overline{x}^i, \overline{u}^i) ds \leq -\xi\varepsilon_i \int_{t_i}^{t_i+\sigma} \alpha_L(\|\overline{x}^i\|) ds < 0, \overline{x}^i \neq 0. \quad (\mathbf{B.2})$$

After time $\delta$ we perform new optimization starting from point $x_{t_i+\delta}$. The value function equals now $V(x^*_{t_i+\delta}) = \overline{T}_{i+1} + \varepsilon_{i+1} \int_{t_i+\delta}^{t_i+\overline{T}_{i+1}+\delta} L(\overline{x}^{i+1}, \overline{u}^{i+1}) ds + \rho_{i+1} q(\overline{x}^{i+1}(t_i + \overline{T}_{i+1} + \delta))$.

Since the pair $\overline{u}^{i+1}$, $\overline{T}_{i+1}$ minimizes $J(u,T;x^*_{t_i+\delta})$ then we obtain inequality

$$V(x^*_{t_i+\delta}) \leq T_{\min} + \varepsilon_{i+1} \int_{t_i+\delta}^{t_i+\overline{T}_i} L(\overline{x}^i, \overline{u}^i) ds + \varepsilon_{i+1} \int_{t_i+\overline{T}_i}^{t_i+T_{\min}+\delta} L(x_s, u_s) ds + \rho_{i+1} q(x_s(t_i + T_{\min} + \delta)).$$

Estimation of the difference in the value function results in

$$V(x^*_{t_i+\delta}) - V(x^*_{t_i}) \leq -(\overline{T}_i - T_{\min}) + (\varepsilon_i + \Delta\varepsilon_i) \int_{t_i+\delta}^{t_i+\overline{T}_i} L(\overline{x}^i, \overline{u}^i) ds + (\varepsilon_i + \Delta\varepsilon_i) \int_{t_i+\overline{T}_i}^{t_i+T_{\min}+\delta} L(x_s, u_s) ds +$$

$$(\rho_i + \Delta\rho_i) q(x_s(t_i + T_{\min} + \delta)) - \varepsilon_i \int_{t_i}^{t_i+\overline{T}_i} L(\overline{x}^i, \overline{u}^i) ds - \rho_i q(\overline{x}^i(t_i + \overline{T}_i)) =$$

$$= -(\overline{T}_i - T_{\min}) - \varepsilon_i \int_{t_i}^{t_i+\delta} L(\overline{x}^i, \overline{u}^i) ds + \Delta\varepsilon_i \left( \int_{t_i+\delta}^{t_i+\overline{T}_i} L(\overline{x}^i, \overline{u}^i) ds + \int_{t_i+\overline{T}_i}^{t_i+T_{\min}+\delta} L(x_s, u_s) ds \right) +$$

$$+ \varepsilon_i \int_{t_i+\overline{T}_i}^{t_i+T_{\min}+\delta} L(x_s, u_s) ds + \rho_i \bigl( q(x_s(t_i + T_{\min} + \delta)) - q(\overline{x}^i(t_i + \overline{T}_i)) \bigr) + \Delta\rho_i q(x_s(t_i + T_{\min} + \delta)).$$

From assumption 1 of the theorem and the fact that $\varepsilon_i \in [0,1]$, $\rho_i \geq 1$ the following inequality occurs $\varepsilon_i \int_{t_i+\overline{T}_i}^{t_i+T_{\min}+\delta} L(x_s, u_s) ds + \rho_i \bigl( q(x_s(t_i + T_{\min} + \delta)) - q(\overline{x}^i(t_i + \overline{T}_i)) \bigr) \leq 0$. Consequently we obtain

$$V(x^*_{t_i+\delta}) - V(x^*_{t_i}) \leq -(\overline{T}_i - T_{\min}) - \varepsilon_i \int_{t_i}^{t_i+\delta} L(\overline{x}^i, \overline{u}^i) ds + \Delta\varepsilon_i \left( \int_{t_i+\delta}^{t_i+\overline{T}_i} L(\overline{x}^i, \overline{u}^i) ds + \int_{t_i+\overline{T}_i}^{t_i+T_{\min}+\delta} L(x_s, u_s) ds \right) +$$

$+ \Delta\rho_i q(x_s(t_i + T_{\min} + \delta))$. From assumption (3.3) we have

$$V(x^*_{t_i+\delta}) - V(x^*_{t_i}) \leq -\xi\varepsilon_i \int_{t_i}^{t_i+\delta} L(\overline{x}^i, \overline{u}^i) ds \leq -\xi\varepsilon_i \int_{t_i}^{t_i+\delta} \alpha_L(\|\overline{x}^i\|) ds < 0, \overline{x}^i \neq 0. \quad (\mathbf{B.3})$$

Combining inequalities **A.1-2**, **B.1-3** will give us

$$V(x^*_{i\delta+\sigma}) + \xi \sum_{k=i_0+1}^{i} \left( \varepsilon_{k-1} \int_{(k-1)\delta}^{k\delta} \alpha_L(\|\overline{x}^{k-1}\|) ds \right) + \varepsilon_i \int_{i\delta}^{i\delta+\sigma} \alpha_L(\|\overline{x}^i\|) ds \leq V(x^*_{i_0\delta}), \quad i = i_0, i_0+1, i_0+2,\ldots.$$

Since we receive trajectory $x^*$ and control $u^*$ in a closed loop system through the concatenation of trajectory $\overline{x}^i$ and controls $\overline{u}^i$ applied to the plant in the time intervals $[t_i, t_{i+1})$ thus we can replace $\overline{x}^i$ with $x^*$

$$V(x^*_{i\delta+\sigma}) + \xi \sum_{k=i_0+1}^{i} \left( \varepsilon_{k-1} \int_{(k-1)\delta}^{k\delta} \alpha_L(\|x^*\|) ds \right) + \varepsilon_i \int_{i\delta}^{i\delta+\sigma} \alpha_L(\|x^*\|) ds \leq V(x^*_{i_0\delta}).$$

Since $0 < \varepsilon_{\min} \leq \varepsilon_i \leq 1$ and $\alpha_L \geq 0$, then

$$V(x_t^*) + \xi\varepsilon_{\min}\int_{i_0\delta}^{t}\alpha_L(\|x^*\|)ds \leq V(x_{i_0\delta}^*) \text{ for } t \geq i_0\delta. \tag{B.4}$$

Therefore $V(x^*(t))$ is limited and optimal horizon sequence is limited i.e. $\overline{T}_i \leq V(x_{i_0\delta}^*)$. Each trajectory $\overline{x}^i$ fulfils the initial condition $\overline{x}^i(t_i) = x^*(t_i)$ and the final condition $\overline{x}^i(t_i + \overline{T}_i) \in \Omega$. From theorem 1.1 we have $\|\overline{x}^i(t)\| \leq (M_1 + M_2\overline{T}_i)\exp(L\overline{T}_i)$ for $t \in [t_i, t_i + \overline{T}_i)$. Hence $\|x^*(t)\| \leq (M_1 + M_2\overline{T}_i)\exp(L\overline{T}_i)$ for $t \in [t_i, t_i + \delta)$ and the consequent $x^*$ is bounded. Since $f$ is continuous then derivative $\dot{x}^* = f(x^*, u^*)$ is also bounded. Inequality (B.4) results in $\lim_{t\to\infty}\int_0^t \alpha_L(\|x^*(s)\|)ds < +\infty$. The following extension of Barbalat's lemma guarantees convergence of the trajectory to zero.

**Lemma 3.1.** Let function $\alpha_L \in \mathcal{K}$ and $x^*(t)$ be an absolutely continuous function on $R_0^+$. If $\|x^*(\cdot)\|_\infty < \infty$, $\|\dot{x}^*(\cdot)\|_\infty < \infty$, $\lim_{T\to\infty}\int_0^T \alpha_L(\|x^*(t)\|)dt < \infty$, to $x^*(t) \to 0$, as $t \to \infty$.

For proof of the lemma see [8, 10].

We will now prove several technical lemmas which in turn will allow us to construct a QTH–RHC algorithm that combines advantages of the time optimal and stabilizing control. The proof of the theorem holds that the value function decreases starting from a certain moment $i_0\delta$ as $\varepsilon_{\min} \leq \varepsilon_i \leq 1$. The lemma below describes behaviour of the value function for $\varepsilon_i \equiv 0$.

**Lemma 3.2.** Let the former assumptions be fulfilled. If $\varepsilon_i \equiv 0$, $\rho = const$, $i = 0,1,2,...$, $\sigma \in [0,\delta]$, then $V(x_{t_i+\sigma}^*) \leq V(x_{t_i}^*) - \min(\sigma, \overline{T}_i - T_{\min})$.

*Proof*: Value function equals $V(x_{t_i}^*) = \overline{T}_i + \rho q(\overline{x}^i(t_i + \overline{T}_i))$. Let $\overline{u}^\sigma$, $\overline{T}_\sigma$ denote a solution to $P(t_i + \sigma, x_{t_i+\sigma}^*)$. Assume that $\overline{T}_i - T_{\min} \geq \delta$. Since the pair $\overline{u}^\sigma$, $\overline{T}_\sigma$ minimizes $J(u;T;x_{t_i+\sigma}^*)$, then limiting the control $\overline{u}^i$ to the interval $[t_i + \sigma, t_i + \overline{T}_i]$ gives us
$V(x_{t_i+\sigma}^*) = \overline{T}_\sigma + \rho q(\overline{x}^\sigma(t_i + \overline{T}_\sigma + \sigma)) \leq \overline{T}_i - \sigma + \rho q(\overline{x}^i(t_i + \overline{T}_i)) = V(x_{t_i}^*) - \sigma$. If $\overline{T}_i - T_{\min} < \delta$, then limiting the control $\overline{u}^i$ to the interval $[t_i + (\overline{T}_i - T_{\min}), t_i + \overline{T}_i]$ gives us
$V(x_{t_i+\sigma}^*) = \overline{T}_\sigma + \rho q(\overline{x}^\sigma(t_i + \overline{T}_\sigma + \sigma)) \leq \overline{T}_i - (\overline{T}_i - T_{\min}) + \rho q(\overline{x}^i(t_i + \overline{T}_i)) = V(x_{t_i}^*) - (\overline{T}_i - T_{\min})$. □

From computational point of view assumption $x(\overline{T}_i + t_i) \in \Omega$ is troublesome to be realized in a calculating algorithm. We will show that it can be fulfilled when choosing big enough penalty coefficient $\rho$.

**Lemma 3.3** Let all former assumptions hold. For any number $\alpha > 0$ there exists a number $\rho > 0$ such that $\|\overline{x}^i(t_i + \overline{T}_i)\|_H^2 \leq \alpha$.

*Proof*: Since the system (1.1) is asymptotically controllable (see def.1.2) then for any number $\eta > 0$ there exists admissible control $u^D : [t_i, t_i + T_D] \to U$, $\overline{T}_i < T_D < +\infty$ such that $\|x^D(t_i + T_D)\| < \eta$ where $x^D$ is the solution of equality (1.1) with the initial condition $x_{t_i}^*$. From the assumption there exists constant $c > 0$ such that $q(x) \geq c\|x\|_H^2$. We obtain the following inequality $c\rho\|\overline{x}^i(t_i + \overline{T}_i)\|_H^2 \leq \rho q(\overline{x}^i(t_i + \overline{T}_i)) \leq \overline{T}_i + \varepsilon_i \int_{t_i}^{t_i+\overline{T}_i} L(\overline{x}^i, \overline{u}^i)ds + \rho q(\overline{x}^i(t_i + \overline{T}_i)) \leq$

$$\leq T_D + \varepsilon_i \int_{t_i}^{t_i+T_D} L(x^D, u^D)ds + \rho q(x^D(t_i + T_D)) \leq T_D + \varepsilon_i \int_{t_i}^{t_i+T_D} L(x^D, u^D)ds + \sup_{\|z\|\leq \eta} \rho q(z).$$ Let us divide both sides by $\rho c$ and we receive $\|\bar{x}^i(t_i + \overline{T}_i)\|_H^2 \leq c^{-1}(\rho^{-1}T_D + \rho^{-1}\varepsilon_i \int_{t_i}^{t_i+T_D} L(x^D, u^D)ds + \sup_{\|z\|\leq \eta} q(z))$.

Taking $\rho > 0$ large enough and accordingly small $\eta > 0$ assures that $\|\bar{x}^i(t_i + \overline{T}_i)\|_H^2 \leq \alpha$. □

**Lemma 3.4** If $\varepsilon_i \equiv 0$, $\rho = const$, $i = 0,1,2,..$, and for certain $j \geq 0$ condition $\dfrac{V(x_{t_j}^*)}{\rho c} \leq \alpha$ is fulfilled, then for every $i \geq j$ there is $\|\bar{x}^i(t_i + \overline{T}_i)\|_H^2 \leq \alpha$.

*Proof*: Value function equals $V(x_{t_i}^*) = \overline{T}_i + \rho q(\bar{x}^i(t_i + \overline{T}_i))$. Evident inequality $\rho q(\bar{x}^i(t_i + \overline{T}_i)) \leq V(x_{t_i}^*)$ occurs. From the assumption there exists constant $c > 0$ such that $q(x) \geq c\|x\|_H^2$. Hence $\|\bar{x}^i(t_i + \overline{T}_i)\|_H^2 \leq c^{-1} q(\bar{x}^i(t_i + \overline{T}_i)) \leq \dfrac{V(x_{t_i}^*)}{\rho c} \leq \alpha$. Lemma **3.2** indicates that the value function does not increase. Hence we obtain inequality

$$\|\bar{x}^{i+1}(t_i + \delta + \overline{T}_{i+1})\|_H^2 \leq \dfrac{V(x_{t_i+\delta}^*)}{\rho c} \leq \dfrac{V(x_{t_i}^*)}{\rho c} \leq \alpha$$ from which the proposition results. □

**Lemma 3.5** If $\varepsilon_{\min} \leq \varepsilon_i \leq 1$, $\rho = const$, $i = 0,1,2,..$, and for certain $j \geq 0$ condition $\dfrac{V(x_{t_j}^*)}{\rho c} \leq \alpha$ is fulfilled then for every $i \geq j$ there is $\|\bar{x}^i(t_i + \overline{T}_i)\|_H^2 \leq \alpha$.

*Proof*: The value function equals $V(x_{t_i}^*) = \overline{T}_i + \varepsilon_i \int_{t_i}^{t_i+\overline{T}_i} L(\bar{x}^i, \bar{u}^i)ds + \rho q(\bar{x}^i(t_i + \overline{T}_i))$. Evident inequality occurs $\rho q(\bar{x}^i(t_i + \overline{T}_i)) \leq V(x_{t_i}^*)$. From the assumption there exists a constant $c > 0$ such that $q(x) \geq c\|x\|_H^2$. Then $\|\bar{x}^i(t_i + \overline{T}_i)\|_H^2 \leq c^{-1} q(\bar{x}^i(t_i + \overline{T}_i)) \leq \dfrac{V(x_{t_i}^*)}{\rho c} \leq \alpha$. From the proof of theorem **3.1** we know that the value function does not increase. Hence we obtain inequality

$$\|\bar{x}^{i+1}(t_i + \delta + \overline{T}_{i+1})\|_H^2 \leq \dfrac{V(x_{t_i+\delta}^*)}{\rho c} \leq \dfrac{V(x_{t_i}^*)}{\rho c} \leq \alpha$$ from which the proposition results. □

For $\varepsilon_i \equiv 0$, the fulfilment of condition $\dfrac{V(x_{t_j}^*)}{\rho c} \leq \alpha$ implies that $x^i(\overline{T}_i + t_i) \in \Omega$ for $i \geq j$. It does not imply however that $x^i(\overline{T}_i + t_i) \in \Omega$ when $\varepsilon_{\min} \leq \varepsilon_i \leq 1$. Lemmas **3.4** and **3.5** indicate that condition $x^i(\overline{T}_i + t_i) \in \Omega$ can be satisfied when choosing a big enough penalty coefficient $\rho$ in the first iteration. Modifications of the penalty coefficient can be necessary in one case only when the coefficient $\varepsilon_i$ gains positive value for the first time.

**4. A Quasi Time Optimal Receding Horizon Control (QTO–RHC).** Theorem **3.1** allows proposition of various schemas stabilizing the system (1.1). The proposed algorithm generates control neighbouring the time optimal control when the state of the system is outside a certain **B** set. When the state of the system reaches this set, the algorithm gradually 'inserts' an integral term into the cost function. This mechanism allows to move in a regular manner from the control neighbouring the time optimal control to the stabilizing control with integral cost function. Let $\mathbf{B} \subset R^n$ be a simply connected point, closed and bounded set including zero. Consider sequential minimization of the cost function (2.1) as $\varepsilon_i \equiv 0$. We have the following

**Lemma 4.1.** If in schema 2.1 $\varepsilon_i \equiv 0$, $\rho = const$ and $\overline{T}_i \geq T_{\min}$, $i = 0,1,2,..$, then the following alternative occurs: Either $\overline{T}_j \in [T_{\min}, T_{\min} + \delta)$ for certain $j \geq 0$ or $\overline{T}_j \geq T_{\min} + \delta$ and $\overline{x}^j(t_j) \in \mathbf{B}$.

*Proof*: Lemma **3.2** indicates that the value function decreases $V(x^*_{t_i+\delta}) \leq V(x^*_{t_i}) - \min(\delta, \overline{T}_i - T_{\min})$. In $k$-th step we have

$$\overline{T}_k + \rho q(\overline{x}(t_k + \overline{T}_k)) = V(x^*_{k\delta}) \leq V(x_0) - \sum_{i=0}^{k-1} \min(\delta, \overline{T}_i - T_{\min}).$$ Since $\overline{T}_i \leq V(x^*_{t_i})$ then from theorem 1.1 we have $\|\overline{x}^i(s)\| \leq (M_1 + M_2 \overline{T}_i) e^{L\overline{T}_i}$. At a certain moment $j\delta$ either the time horizon $\overline{T}_j$ will decrease so that $\overline{T}_j \in [T_{\min}, T_{\min} + \delta)$ or the set $\mathbf{B}$ will be reached. □

Speaking inaccurately, it is either not possible to decrease the horizon of $\delta$ or set $\mathbf{B}$ has been reached. Utilizing all former results it is possible to propose the following algorithm.

**Algorithm 4.1.** Given: Cost function (2.1), plant model (1.1), sets $\mathbf{B}$ and $\Omega$, number $\alpha > 0$ such that the following implication occurs $c\|x\|_H^2 \leq \alpha \Rightarrow x \in \Omega$, $\xi \in (0,1)$ (e.g. $\xi = 0.1$), $\rho \geq 1$ (e.g. $\rho = 100$). Substitute $i = 0$, $t_i = 0$, $\varepsilon_i = 0$.

**1.** Make estimation of the state $x_{t_i}$;

**2.** Solve $P(t_i, x_{t_i})$;

**2a**. If $\dfrac{V(x_{t_i})}{\rho c} > \alpha$ then substitute $\rho \leftarrow \gamma \rho$, $\gamma > 1$ (e.g. $\gamma = 2$) and go to **2**;

**2b**. If $\overline{T}_i = T_{\min} \wedge \varepsilon_i = 0$ then choose $\varepsilon_i \in (0,1)$ (e.g. $\varepsilon_i = 0.01$) and go to **2**;

**3.** Apply the initial part of control $\overline{u}^i$ in time interval $t \in [t_i, t_i + \delta)$;

**4**. Make estimation of the state of $x_{t_i}$;

**5**. Update coefficient $\varepsilon_i$ according to the formula $\varepsilon_{i+1} = \begin{cases} \min(\varepsilon_i + \Delta\overline{\varepsilon}_i, 1), & \text{if } x_{t_i} \in \mathbf{B} \\ \varepsilon_i, & \text{if } x_{t_i} \notin \mathbf{B} \end{cases}$ (4.1)

where

$$\Delta\overline{\varepsilon}_i = (1-\xi) \frac{\delta + \varepsilon_i \int_{t_i}^{t_i+\delta} L(\overline{x}^i, \overline{u}^i) dt}{\int_{t_i+\delta}^{t_i+\overline{T}_i} L(\overline{x}^i, \overline{u}^i) dt} \quad \text{for } \overline{T}_i \geq T_{\min} + \delta,$$ (4.2)

and

$$\Delta\overline{\varepsilon}_i = (1-\xi) \frac{(\overline{T}_i - T_{\min}) + \varepsilon_i \int_{t_i}^{t_i+\delta} L(\overline{x}^i, \overline{u}^i) ds}{\int_{t_i+\delta}^{t_i+\overline{T}_i} L(\overline{x}^i, \overline{u}^i) ds + q(\overline{x}^i(t_i + \overline{T}_i))} \quad \text{for } \overline{T}_i \in [T_{\min}, T_{\min} + \delta);$$ (4.3)

**6.** Replace $t_{i+1} \leftarrow t_i + \delta$, $i \leftarrow i+1$ and go to **2**.

**Theorem 4.1.** Let assumptions 1 and 2 of theorem **3.1** and all former ones hold. Then trajectory $x^*(t)$ generated by algorithm **4.1** converges asymptotically to zero.

*Proof*: We will show that assumptions of theorem **3.1** are fulfilled. Increase the penalty coefficient $\rho$ at point 2 and 2a of the algorithm until condition $\dfrac{V(x_{t_j})}{\rho c} \le \alpha$ is fulfilled. Lemma 3.3 indicates that this procedure must finish at a certain moment and optimal trajectory will then satisfy the condition $\bar{x}^i(\bar{T}_i + t_i) \in \Omega$. Hence, control is admissible in the sense of definition 1.1. The condition at point 2b guarantees fulfilment of assumption 4 of theorem 3.1 in case when ($\bar{T}_i = T_{\min} \wedge \varepsilon_i = 0$). The formula for updating the coefficient $\varepsilon_i$ in point 5 and lemma 4.1 guarantee fulfilment of assumption 4 of the theorem in the remaining cases. Lemmas 3.4 and 3.5 indicate that a penalty coefficient $\rho$ does not alter except when ($T_i = T_{\min} \wedge \varepsilon_i = 0$). In this exceptional case the penalty coefficient may (but does not have to) increase. However, it does not change the fact that from a certain moment this coefficient remains constant. From equality 4.2 we immediately obtain

$$\Delta \varepsilon_i \int_{t_i+\delta}^{t_i+\bar{T}_i} L(\bar{x}^i, \bar{u}^i) ds \le \delta + (1-\xi)\varepsilon_i \int_{t_i}^{t_i+\delta} L(\bar{x}^i, \bar{u}^i) ds.$$

Since $\Delta \rho_i = 0$ then (3.2) occurs. Assumption 1 of theorem 3.1 results in the following inequality $\int_{t_i+\bar{T}_i}^{t_i+T_{\min}+\delta} L(x_s, u_s) ds \le q(\bar{x}^i(t_i + \bar{T}_i))$. From the above inequality and inequality (4.3) we obtain

$$\Delta \varepsilon_i = (1-\xi) \dfrac{(\bar{T}_i - T_{\min}) + \varepsilon_i \int_{t_i}^{t_i+\delta} L(\bar{x}^i, \bar{u}^i) ds}{\int_{t_i+\delta}^{t_i+\bar{T}_i} L(\bar{x}^i, \bar{u}^i) ds + q(\bar{x}^i(t_i + \bar{T}_i))} \le \dfrac{(\bar{T}_i - T_{\min}) + (1-\xi)\varepsilon_i \int_{t_i}^{t_i+\delta} L(\bar{x}^i, \bar{u}^i) ds}{\int_{t_i+\delta}^{t_i+\bar{T}_i} L(\bar{x}^i, \bar{u}^i) ds + \int_{t_i+\bar{T}_i}^{t_i+T_{\min}+\delta} L(x_s, u_s) ds}.$$

Since $\Delta \rho_i = 0$ then (3.3) occurs and this finishes the proof. □

Coefficient $\varepsilon_i$ equals zero in the initial phase of the algorithm. Increase of $\varepsilon_i$ begins at the moment of reaching set **B** or in case when $T_i = T_{\min} \wedge \varepsilon_i = 0$.

**5. Stabilizable systems with quadratic cost function.** In general, it is not an easy task to define a set $\Omega$ and a function $q$ that fulfil the assumptions of theorem 3.1. However, in case when the integrand in (2.1) is quadratic and with the assumption of stabilizabilty of the linear system in the neighbourhood of zero it is possible to give a general algorithm [2,3,1]. Assume that $L(x,u) = x^T W x + u^T R u$, $W^T = W > 0$, $R^T = R > 0$, $q(x) = k x^T H x$, $k > 1$. Let $A = f_x^T(0,0)$, $B = f_u^T(0,0)$. Also assume that the pair $(A,B)$ is stabilizable and the pair $(W,A)$ is detectable. Then there exists such matrix $K$ that matrix $A_K = A + BK$ is exponentially stable. Moreover, Lyapunov equality $A_K^T H + H A_K = -(W + K^T R K)$ has a unique solution $H = H^T > 0$. Let the stabilizing control $u_s$ (see theorem 3.1) be $u_s = K x_s$. The following lemma shows how to fulfil assumption 1 of theorem 3.1.

**Lemma 5.1.** Let $\Omega = \{x \in R^n ; x^T H x \le \alpha\}$ and function $x_s : R_0^+ \to R^n$ be the solution of equality (1.1) with an initial condition $x_s(0) \in \Omega$ and control $u_s = K x_s$.
Then there exists such a number $\alpha > 0$ that:
**1.** Inequality $\dot{q}(x_s(\tau)) + L(x_s(\tau), u_s(\tau)) \le 0$ is fulfilled in set $\Omega$.
**2.** $\Omega$ is the region of attraction of the equilibrium point. Each trajectory of the system (1.1) with control $u_s$ starting from $\Omega$ remains in $\Omega$.
**3.** Control satisfies constraints (1.2) i.e. $K x_s(\tau) \in U$.

The way to determine the number $\alpha$ and the proof of the lemma can be fund in works [1,2]. Set $\Omega$ is an *n*–dimensional ellipse with the centre in point 0. Let us test assumptions for problem $P(t_i, x_{t_i})$. Assuming that $\alpha_L(\|x\|) = \lambda_{\min}(W)\|x\|^2$ shows that $L(x,u) \geq \alpha_L(\|x\|)$ and that $\alpha_L \in \mathcal{K}$. Since $k > 1$ then $q(x) \geq c\|x\|_H^2$ with a constant $c = 1$. The condition $c\|x\|_H^2 \leq \alpha$ evidently implies that $x \in \Omega$.

**6. Examples.** We will show application of the algorithm 4.1 to control non linear systems. In all examples the cost function was quadratic. Set $\Omega$ and function $q$ were determined by the method described in chapter 5. A version of Monotone Structural Evolution [11] algorithm modified by the author was used to solve the optimal control problems.

**Example 6.1.** Consider the control of pendulum described with equalities $\dot{x}_1 = x_2$, $\dot{x}_2 = \sin(x_1) + 0.3u$, where $x_1$ - pendulum angle, $x_2$ - pendulum speed. The task is to guide the pendulum from a lower stable equilibrium point $x_0 = [-\pi\ 0]^T$ to an upper unstable equilibrium point $x_r = [0\ 0]^T$. The upper and lower values of control were $u_{\max} = 1$, $u_{\min} = -1$. The stabilizing controller was determined by solving LQ problem for system linearized in point $x_r$, with matrices $R = 500$, $W = \text{diag}[500\ 500]$. Gain matrix of the controller was equal to $K = [-6.81\ -6.81]$. Set $\Omega = \{x \in R^2; x^T H x \leq 0.01\}$, $\mathbf{B} = \{x \in R^2; |x_i| \leq 0.3, i = 1,2\}$, $\delta = 0.05$, $T_{\min} = 0.5$. The first numerical experiment solved the time optimal problem. The time optimal horizon equalled $T_D = 12.56$. The time optimal control and its corresponding trajectory were compared to control $u^*(t)$ and trajectory $x^*(t)$ generated by algorithm 4.1. The results of the simulation are presented in fig. 6.1-8.

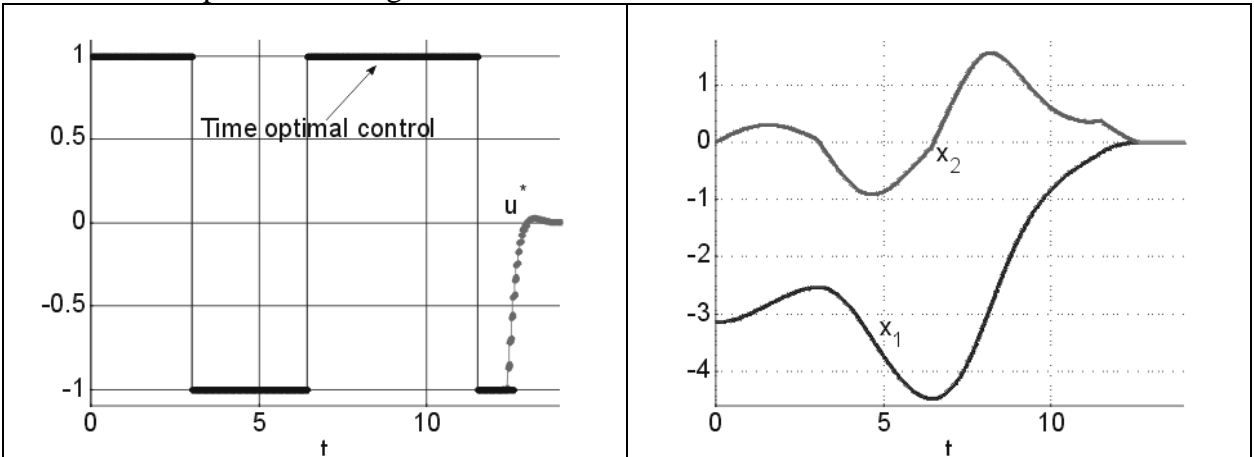

Fig. 6.1.1. Control $u^*(t)$ and time optimal control.　　Fig. 6.1.2 Time optimal trajectory and trajectory $x^*(t)$.

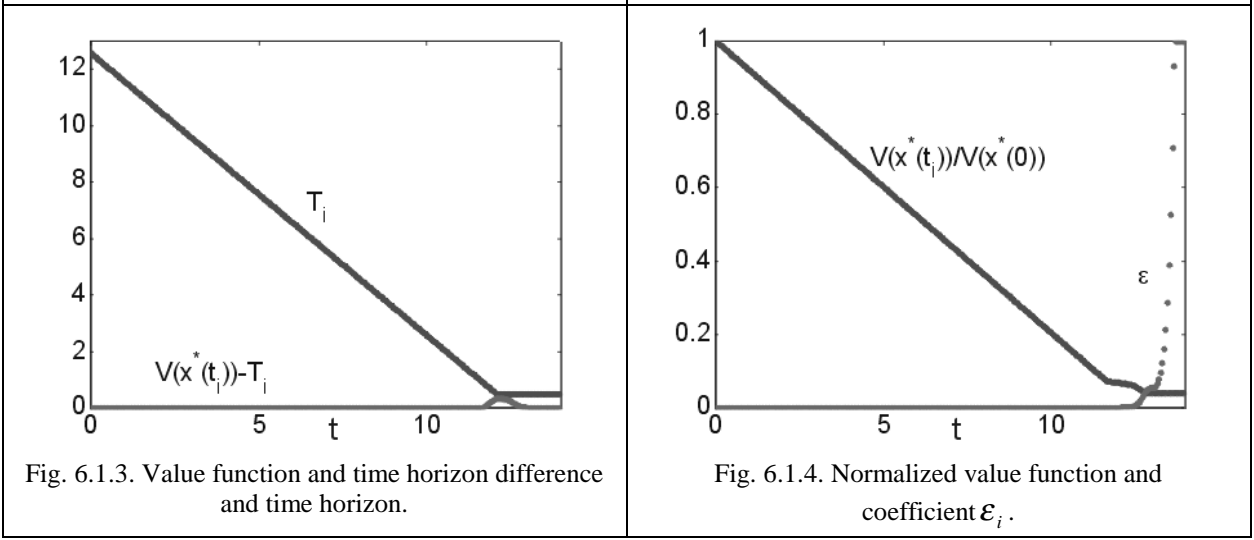

Fig. 6.1.3. Value function and time horizon difference and time horizon.　　Fig. 6.1.4. Normalized value function and coefficient $\varepsilon_i$.

In order to test the algorithm more precisely in the neighbourhood of the equilibrium point, figs. 6.1.5-8. show the final part of trajectory and control.

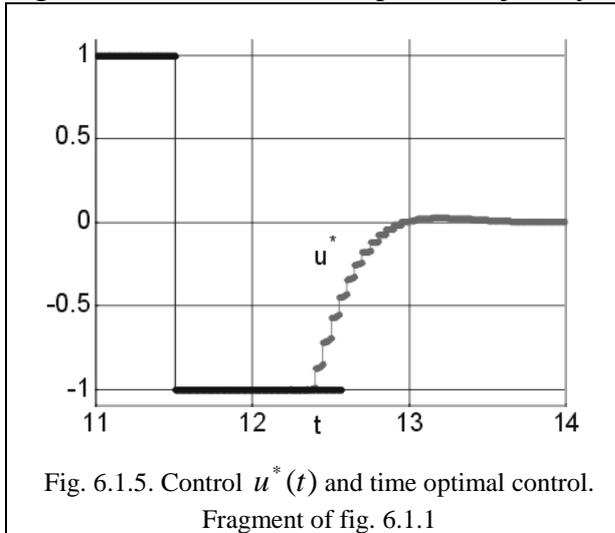

Fig. 6.1.5. Control $u^*(t)$ and time optimal control. Fragment of fig. 6.1.1

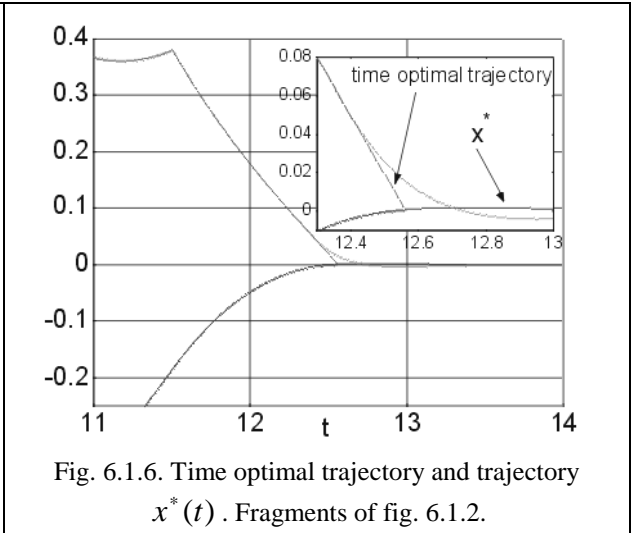

Fig. 6.1.6. Time optimal trajectory and trajectory $x^*(t)$. Fragments of fig. 6.1.2.

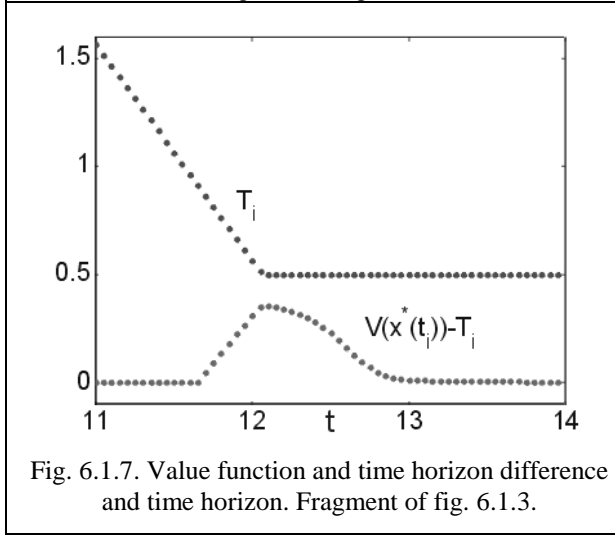

Fig. 6.1.7. Value function and time horizon difference and time horizon. Fragment of fig. 6.1.3.

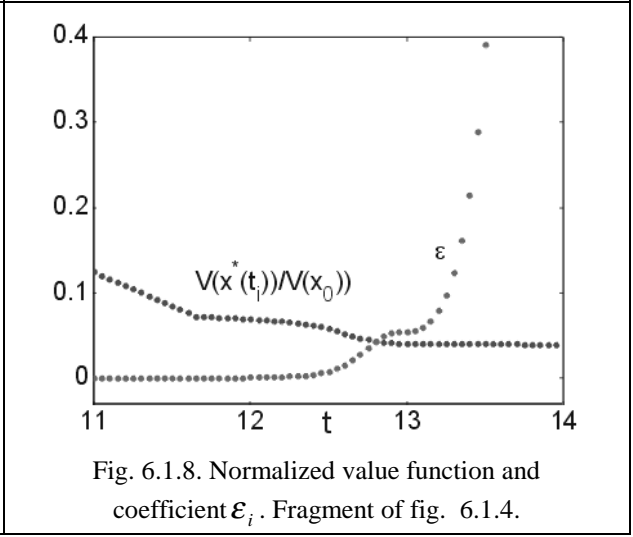

Fig. 6.1.8. Normalized value function and coefficient $\varepsilon_i$. Fragment of fig. 6.1.4.

**Example 6.2.** Stabilization of a pendulum on a cart system. The system of a pendulum on a cart is a known and well examined example of a nonlinear system [12,13]. In Turnau work [12] the following non dimensional equations were introduced: $\dot{x}_1 = x_3$, $\dot{x}_2 = x_4$,

$$\dot{x}_3 = \frac{v_1(x,u) + v_2(x)\cos x_2}{d(x)}, \quad \dot{x}_4 = \frac{v_1(x,u)\cos x_2 + c_4 v_2(x)}{d(x)}, \quad d(x) = c_4 - \cos^2 x_2,$$

$v_1(x,u) = u - x_4^2 \sin x_2 - b_2 x_3$, $v_2(x) = \sin x_2 - b_3 x_4$, where $x_1$ – cart position, $x_2$ – pendulum angle, $x_3$ – cart speed, $x_4$ – pendulum speed. Parameters of the model were as follows: $b_2 = 0.4256$, $b_3 = 3.1564 \cdot 10^{-4}$, $c_4 = 11.2135$, $u_{max} = 3.9351$, $u_{min} = -3.9351$. The goal of the control was to stabilize the pendulum in an upper unstable equilibrium point $x_r = [0,0,0,0]$. The stabilizing controller was determined by solving LQ problem for system linearized in point $x_r$, with matrices $R = 6.46$, $W = \text{diag}[1132, 100, 1, 1]$. Gain matrix of the controller equalled $K = [13.24, -81.74, 43.65, -80.63]$. The initial condition equalled $x_0 = [0\ 0.04\ 0\ 0]$, $\Omega = \{x \in R^4; x^T H x \leq 0.07\}$, $\mathbf{B} = \{x \in R^4; |x_i| \leq 0.05, i = 1,...,4\}$, $\delta = 0.2$, $T_{min} = 1.3$, $\varepsilon_0 = 0.05$. The results of the simulation are presented in figs. 6.2.1-4.

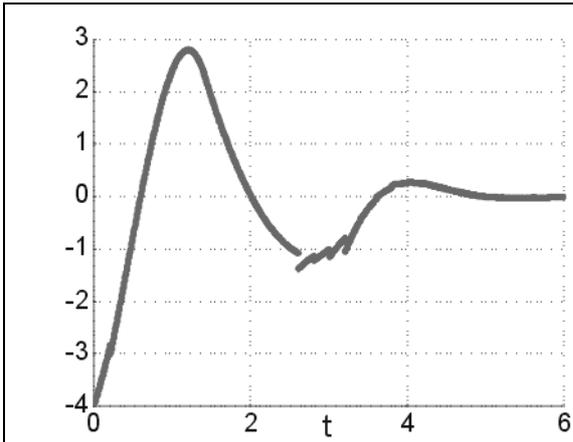

Fig. 6.2.1. Control.

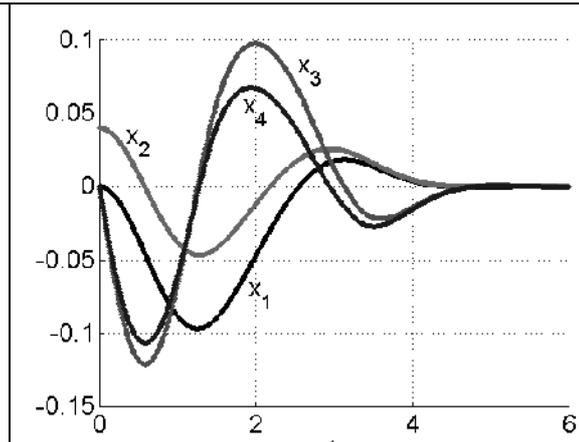

Fig. 6.2.2. Trajectory.

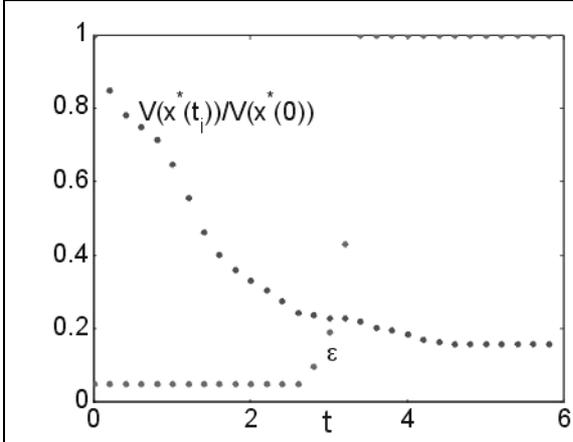

Fig. 6.2.3. Normalized value function and coefficient $\varepsilon_i$.

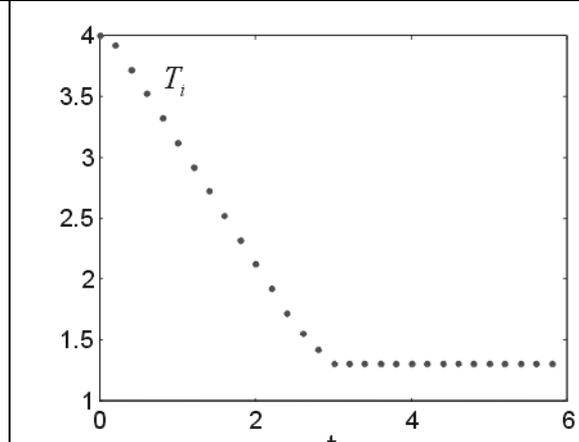

Fig. 6.2.4. Time horizon.

The second experiment compared the performance of a receding horizon algorithm to that of the LQ controller with an infinite control time. All parameters were the same as in the previous experiment. The results of the simulation are presented in figs. 6.2.5-6.

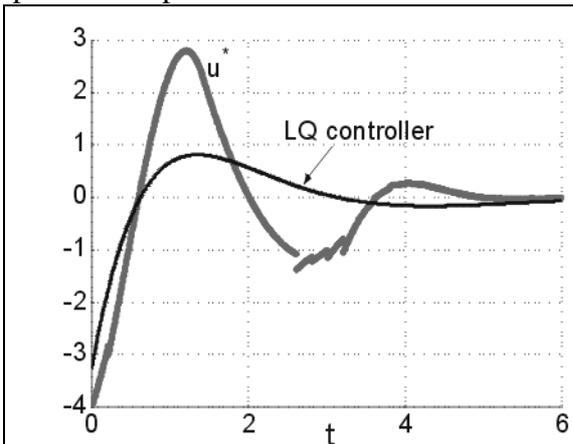

Fig. 6.2.5. Comparison of control attained through the performance of algorithm 4.1 to control generated by LQ algorithm.

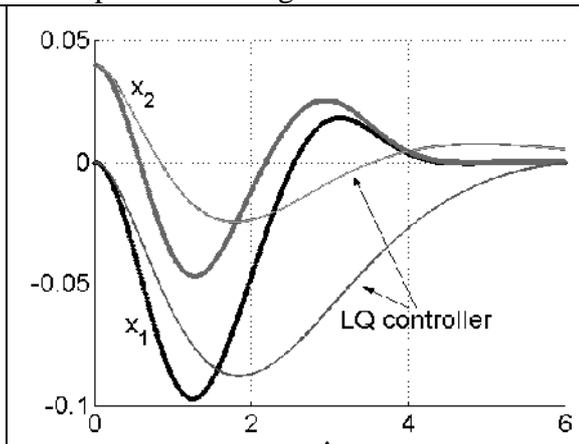

Fig. 6.2.6. Comparison of trajectory received as a result of the performance of algorithm 4.1 to trajectory generated by LQ algorithm. For figure legibility only the two first components of the state are shown.

**Conclusions.** The paper presents a Quasi Time Optimal Receding Horizon Control (QTO–RHC) algorithm. The proposed algorithm generates control neighbouring to time optimal control when the state of the system is far from the equilibrium point. During the control process the state of the system approaches a certain set **B** defined by the user. When the state of the system reaches set **B**, the adaptation of cost function begins. The goal of this adaptation is to move from a time optimal control to a stabilizing control. If the cost function is quadratic then the algorithm converges asymptotically to the LQ controller. Sufficient conditions for stability of the closed loop system have been determined and the manner of the adaptation of the cost function has been described. The simulating experiments show that the algorithm allows for a very regular transition from time optimal trajectory to the stabilization phase (see fig. 6.1.5-6). Since the time horizon decreases until it reaches the minimal value hence the calculating complexity when solving the optimal control problem also decreases. Analysis of the simulation results indicates that the growth of the coefficient $\varepsilon_i$ is relatively fast thereby the adaptation of the cost function occurs surprisingly fast. Comparison of the performance of the algorithm 4.1 to the performance of the LQ controller indicates significant shortening of control time (see fig. 6.2.10). The advantage of the proposed approach is the fact that the user can arbitrary define set **B** and the minimal horizon. Analysis of the algorithm robustness should be the goal of further research. It seems that Ledayev [7] proposition offers promising perspectives.

**Acknowledgments**: This paper was supported by the USC research grant No. 10120038.

**References.**

[1] Bania P. 2005: *Determination of terminal constraints in nonlinear predictive control*. Proc. of fifth national conference Computer Methods and Systems CMS'05, Kraków, Poland. In Polish

[2] Chen H., Allgöwer F. 1998: *A Quasi Infinite Horizon Nonlinear Model Predictive Control Scheme With Guaranted Stability* Automatica vol. 34, No. 10, pp. 1205-1217.

[3] Findeisen R., Imsland L., Allgöwer F., Foss B.A., 2003a: *State and Output Feedback Nonlinear Model Predictive Control: An Overview* Europian Journal of Control 9:190-2006.

[4] Findeisen R., Imsland L., Allgöwer F., Foss B.A., 2003b: *Output Feedback stabilization of Consrtained Systems With Nonlinear Predictive Control* International Journal of Robust and Nonlinear Control 2003;13:211-227.

[5] Fontes A. C. C., Magni L. 2003: *Min–Max Model Predictive Control of Nonlinear Systems Using Discontinuous Feedbacks* IEEE Transactions on Automatic Control vol. 48, NO. 10, October 2003 pp. 1750 –1755.

[6] Fontes F. A., 2000: *A general Framework to Design Stabilizing nonlinear Model Predictive Controller* Syst. Contr. Letters 42(2):127-143.

[7] Ledayev Y. S. 2002:*Robustnes of dicontinuous feedback in control under disturbance* : in Proc. of the 41[st] IEEE Conference on Decision and Control Las Vegas, Nevada USA, December 2002.

[8] Logemann H. Ryan E. P. 2004: *Asymptotic Behaviour of Nonlinear Systems* in The Mathematical Association of America Monthly December 2004, (http://www.maa.org).

[9] Mayne D. Q., Michalska H. 1990: *Receding Horizon Control of Nonlinear Systems*. IEEE Transactions on Automatic Control, 35:814-824.

[10] Michalska H., Mayne D. Q. 1993: *Robust Receding Horizon Control of Constrained Nonlinear Systems*. IEEE Transactions on Automatic Control, 38:1623-1633.

[11] Szymkat M., Korytowski A. 2003: *Method of monotone structural evolution for control and state constrained optimal control problems* in Proc. of the Europian Control Conference ECC 2003, 1-4 September University of Cambridge, UK.

[12] Turnau A. 2002: *Target Control of nonlinear systems in real time –intelligent and time-optimal algorithms*. AGH Kraków 2002 Poland. In Polish.

[13] Turnau A., Korytowski A., Szymkat M. 1999: *Time optimal control for the pendulum-cart system in real-time*. Proc. 1999 IEEE CCA, Kohala Coast, Hawai'i, August 22-27, 1999, 1249-1254.

[14] http://home.agh.edu.pl/~pba/thm11.pdf